\documentclass[letterpaper, 10 pt, conference]{ieeeconf}
\usepackage{amsmath ,amssymb,times,mathptmx,epsfig,graphicx,color}
\graphicspath{{./},{./epsfigures/}}

\begin{document}
\newtheorem{thm}{Theorem}
\newtheorem{cor}[thm]{Corollary}
\newtheorem{conj}[thm]{Conjecture}
\newtheorem{prop}[thm]{Proposition}
\newtheorem{problem}{Problem}
\newtheorem{remark}[thm]{Remark}
\newtheorem{defn}[thm]{Definition}
\newtheorem{ex}[thm]{Example}

\IEEEoverridecommandlockouts 
\newcommand{\mR}{{\mathbb R}}
\newcommand{\mP}{{\mathbb P}}
\newcommand{\D}{{\mathbb D}}
\newcommand{\E}{{\mathbb E}}
\newcommand{\mcN}{{\mathcal N}}
\newcommand{\mcR}{{\mathcal R}}
\newcommand{\tr}{\operatorname{trace}}
\newcommand{\ignore}[1]{}

\overrideIEEEmargins  
\newcommand{\mike}{\color{magenta}}
\newcommand{\rike}{\color{red}}

\title{\LARGE \bf Stochastic control, entropic interpolation and gradient flows on Wasserstein product spaces*}

\author{Yongxin Chen$^{1}$, Tryphon Georgiou$^{2}$ and Michele Pavon$^{3}$%
\thanks{*Research partially supported by the NSF under Grant ECCS-1509387, the AFOSR under Grants FA9550-12-1-0319, and FA9550-15-1-0045 and by the University of Padova Research Project CPDA 140897.}%
\thanks{$^{1}$Yongxin Chen is with the Department of Mechanical Engineering,
University of Minnesota, Minneapolis, Minnesota MN 55455, USA   
{\tt\small chen2468@umn.edu}}%
\thanks{$^{2}$ Tryphon Georgiou is with the Department of Electrical and Computer Engineering,
University of Minnesota, Minneapolis, Minnesota MN 55455, USA 
{\tt\small tryphon@umn.edu}}%
\thanks{$^{3}$ Michele Pavon is with the Dipartimento di Matematica,
Universit\`a di Padova, via Trieste 63, 35121 Padova, Italy  
{\tt\small pavon@math.unipd.it}}%
}

\maketitle
\thispagestyle{empty}
\pagestyle{empty}

\begin{abstract}
Since the early nineties, it has been observed that the Schr\"{o}dinger bridge problem can be formulated as a stochastic control problem with atypical boundary constraints. This in turn has a fluid dynamic counterpart where the flow of probability densities represents an {\em entropic interpolation} between the given initial and final marginals. In the zero noise limit, such entropic interpolation converges in a suitable sense to the {\em displacement interpolation} of optimal mass transport (OMT).
We consider two absolutely continuous curves in Wasserstein space ${\cal W}_2$ and  study the evolution of the relative entropy on ${\cal W}_2\times {\cal W}_2$ {\em on a finite time interval}. Thus, this study differs from previous work in OMT theory concerning relative entropy from a fixed (often equilibrium) distribution (density). We derive a gradient flow on Wasserstein product space. We find the remarkable property that fluxes in the two components are {\em opposite}. Plugging in the ``steepest descent" into the evolution of the relative entropy we get what appears to be a new formula: The two flows approach each other at a faster rate than that of two solutions of the same Fokker-Planck. We then study the evolution of relative entropy in the case of uncontrolled-controlled diffusions. In two special cases of the Schr\"{o}dinger bridge problem, we show that such relative entropy may be monotonically decreasing or monotonically increasing.
\end{abstract}

\section{Introduction}
In the Schr\"{o}dinger bridge problem (SBP) \cite{F2}, one seeks the random evolution (a probability measure on path-space) which is closest in the relative entropy sense to a prior Markov diffusion evolution and has certain prescribed initial and final marginals $\mu$ and $\nu$. As already observed by Schr\"{o}dinger \cite{S1,S2}, the problem may be reduced to a {\em static} problem which, except for the cost, resembles the Kantorovich relaxed formulation of the optimal mass transport problem (OMT). Considering that since \cite{BB} (OMT) also has a dynamic formulation, we have two problems which admit equivalent static and dynamic versions \cite{leo2}. Moreover, in both cases, the solution entails a flow of one-time marginals joining $\mu$ and $\nu$. The OMT yields a {\em displacement interpolation flow} whereas the SBP provides an {\em entropic interpolation flow}. 

Trough the work of Mikami, Mikami-Thieullen and Leonard \cite{Mik, mt, MT,leo,leo2}, we know that the OMT may be viewed as a ``zero-noise limit" of SBP when the prior is a sort of uniform measure on path space with vanishing variance. This connection has been extended to more general prior evolutions in \cite{CGP3,CGP4}. Moreover, we know that, thanks to a very useful intuition by Otto \cite{O}, the displacement interpolation flow $\{\mu_t; 0\le t\le 1\}$ may be viewed as a constant-speed geodesic joining $\mu$ and $\nu$ in Wasserstein space \cite{Vil}.  
What can be said from this geometric viewpoint of the entropic flow? It cannot be a geodesic, but can it be characterized as a curve minimizing a suitable action? In \cite{CGP3}, we showed that this is indeed the case resorting to a time-symmetric fluid dynamic formulation of  SBP. The action features an extra term which is a {\em Fisher information functional}. Moreover, this characterization of the Schr\"{o}dinger bridge answers at once a question posed by Carlen \cite[pp. 130-131]{Carlen}.

It has been observed since the early nineties that SBP can be turned, thanks to {\em Girsanov's theorem},  into a stochastic control problem with atypical boundary constraints, see \cite{DP,Bl,DPP,PW,FHS}. The latter has a fluid dynamic counterpart. It is therefore interesting to compare the flow associated to the uncontrolled evolution (prior) to the optimal one. In particular, it is interesting to study the evolution of the relative entropy on the {\em product} Wasserstein space  {\em on a finite time interval}. Thus, this study differs from previous work in OMT theory concerning relative entropy from an equilibrium distribution (density). We derive in Section \ref{REWP} a gradient flow on Wasserstein product space. We find the remarkable property that fluxes in the two components are {\em opposite}. Plugging in the``steepest descent" into the evolution of the relative entropy we get what appears to be a new formula (\ref{REFF}): The two flows approach each other at a faster rate than that of two solutions of the same Fokker-Planck. We then study the evolution of relative entropy in the case of uncontrolled-controlled diffusions. We show by one special case of the Schr\"{o}dinger bridge problem that such relative entropy may even be monotonically increasing.

The paper is outlined as follows. In Section \ref{OMT}, we recall some fundamental facts and concepts from the theory of optimal transportation. In Section \ref{FPGR}, we review the variational formulation of the Fokker-Planck equation as a gradient flow on Wasserstein space. Section \ref{REWP}, we study the evolution of relative entropy on Wasserstein product space.  In Section \ref{NF}, we recall some basic elements of the Nelson-F\"{o}llmer kinematics of finite-energy diffusions. Finally, in Section \ref{stochastic control}, we study the relative entropy change in the case of a controlled evolution. This is then specialized to the Schr\"{o}dinger bridge.

\section{Elements of optimal mass transport theory}\label{OMT}
The literature on this problem is by now so vast and our degree of competence is such that we shall not even attempt here to give a reasonable and/or balanced introduction to the various fascinating aspects of this theory. Fortunately, there exist excellent monographs and survey papers on this topic, see \cite{RR,E,Vil,AGS,Vil2, OPV}, to which we refer the reader. We shall only briefly review some concepts and results which are relevant for the topics of this paper.
\subsection{The static problem}

Let $\nu_0$ and $\nu_1$ be probability measures on the measurable spaces $X$ and $Y$, respectively. Let $c:X\times Y\rightarrow [0,+\infty)$ be a measurable map with $c(x,y)$ representing the cost of transporting a unit of mass from location $x$ to location $y$.  Let ${\cal T}_{\nu_0\nu_1}$ be the family of measurable maps $T:X\rightarrow Y$ such that $T\#\nu_0=\nu_1$, namely such that $\nu_1$ is the {\em push-forward} of $\nu_0$ under $T$. Then Monge's optimal mass transport problem (OMT) is
\begin{equation}\label{monge}\inf_{T\in {\cal T}_{\nu_0\nu_1}}\int_{X\times Y} c(x,T(x))d\nu_0(x).
\end{equation}
As is well known, this problem may be unfeasible, namely the family ${\cal T}_{\nu_0\nu_1}$ may be empty. This is never the case for the ``relaxed" version of the problem studied by Kantorovich in the 1940's
\begin{equation}\label{kantorovich}
\inf_{\pi\in\Pi(\nu_0,\nu_1)}\int_{X\times Y} c(x,y)d\pi(x,y)
\end{equation}
where $\Pi(\nu_0,\nu_1)$ are ``couplings" of $\nu_0$ and $\nu_1$, namely probability distributions on $X\times Y$ with marginals $\nu_0$ and $\nu_1$. Indeed, $\Pi(\nu_0,\nu_1)$ always contains the product measure $\nu_0\otimes\nu_1$. Let us specialize the Monge-Kantorovich problem (\ref{kantorovich}) to the case $X=Y=\mR^N$ and $c(x,y)=\|x-y\|^2$. Then, if $\nu_1$ does not give mass to sets of dimension $\le n-1$, by Brenier's theorem \cite[p.66]{Vil}, there exists a unique optimal transport plan $\pi$ (Kantorovich) induced by a $d\nu_0$ a.e. unique map $T$ (Monge), $T=\nabla\varphi$, $\varphi$ convex, and we have 
\begin{equation}\label{optmap}\pi=(I\times\nabla\varphi)\#\nu_0, \quad \nabla\varphi\#\nu_0=\nu_1.
\end{equation}
Here $I$ denotes the identity map. Among the extensions of this result, we mention that to strictly convex, superlinear costs $c$ by Gangbo and McCann \cite{GM}.
The optimal transport problem may be used to introduce a useful distance between probability measures. Indeed, let $\mathcal P_2(\mR^N)$ be the set of probability measures $\mu$ on $\mR^N$ with finite second moment. For $\nu_0, \nu_1\in\mathcal P_2(\mR^N)$, the Wasserstein (Vasershtein) quadratic distance, is defined by
\begin{equation}\label{Wasserdist}
W_2(\nu_0,\nu_1)=\left(\inf_{\pi\in\Pi(\nu_0,\nu_1)}\int_{\mR^N\times\mR^N}\|x-y\|^2d\pi(x,y)\right)^{1/2}.
\end{equation}
As is well known \cite[Theorem 7.3]{Vil}, $W_2$ is a {\em bona fide} distance. Moreover, it provides a most natural way to  ``metrize" weak convergence\footnote{$\mu_k$ converges weakly to $\mu$ if $\int_{\mR^N}fd\mu_k\rightarrow\int_{\mR^N}fd\mu$ for every continuous, bounded function $f$.} in $\mathcal P_2(\mR^N)$ \cite[Theorem 7.12]{Vil}, \cite[Proposition 7.1.5]{AGS} (the same applies to the case $p\ge 1$ replacing $2$ with $p$ everywhere). The {\em Wasserstein space} $\mathcal W_2$ is defined as the metric space $\left(\mathcal P_2(\mR^N),W_2\right)$. It is a {\em Polish space}, namely a separable, complete metric space.

\subsection{The dynamic problem}

So far, we have dealt with {\em the static} optimal transport problem. Nevertheless, in \cite[p.378]{BB} it is observed that ``...a continuum mechanics formulation was already implicitly contained in the original problem addressed by Monge... Eliminating the time variable was just a clever way of reducing the dimension of the problem". Thus, a {\em dynamic} version of the OMT problem was already {\em in fieri} in Gaspar Monge's 1781 {\em ``M\'emoire sur la th\'eorie des d\'eblais et des remblais"}\,! It was elegantly accomplished by Benamou and Brenier in \cite{BB} by showing that 
\begin{subequations}\label{eq:BB}
\begin{eqnarray}\label{BB1}&&W_2^2(\nu_0,\nu_1)=\inf_{(\mu,v)}\int_{0}^{1}\int_{\mR^N}\|v(x,t)\|^2\mu_t(dx)dt,\\&&\frac{\partial \mu}{\partial t}+\nabla\cdot(v\mu)=0,\label{BB2}\\&& \mu_0=\nu_0, \quad \mu_1=\nu_1.\label{boundary}
\end{eqnarray}\end{subequations}
Here the flow $\{\mu_t; 0\le t\le 1\}$ varies over continuous maps from $[0,1]$ to $\mathcal P_2(\mR^N)$ and $v$ over smooth fields.  In \cite{Vil2}, Villani states at the beginning of Chapter $7$ that two main motivations for the time-dependent version of OMT are
\begin{itemize}
\item a time-dependent model gives a more complete description of the
transport;
\item the richer mathematical structure will be useful later on.
\end{itemize}
We can add three further reasons:
\begin{itemize}
\item it opens the way to establish a connection with the {\em Schr\"{o}dinger bridge} problem, where the latter appears as a regularization of the former \cite{Mik, mt, MT,leo,leo2,CGP3,CGP4};
\item it allows to view the optimal transport problem as an (atypical) {\em optimal control} problem \cite{CGP}-\cite{CGP4}.
\item In some applications such as interpolation of images \cite{CGP5} or spectral morphing \cite{JLG}, the interpolating flow is essential!
\end{itemize}
Let $\{\mu^*_t; 0\le t\le 1\}$ and $\{v^*(x,t); (x,t)\in\mR^N\times[0,1]\}$ be optimal for (\ref{eq:BB}). Then 
$$\mu^*_t=\left[(1-t)I+t\nabla\varphi\right]\#\nu_0, $$ with $T=\nabla\varphi$ solving Monge's problem,  provides, in McCann's language, the {\em displacement interpolation} between $\nu_0$ and $\nu_1$. Then $\{\mu^*_t; 0\le t\le 1\}$ may be viewed as a constant-speed geodesic joining $\nu_0$ and $\nu_1$ in Wasserstein space (Otto). This formally endows  $\mathcal W_2$ with a ``pseudo" Riemannian structure.
McCann discovered \cite{McCann} that certain functionals are {\em displacement convex}, namely convex along Wasserstein geodesics. This has led to a variety of applications. Following one of Otto's main discoveries \cite{JKO,O}, it turns out that a large class of PDE's may be viewed as {\em gradient flows}  on the Wasserstein space ${\cal W}_2$. This interpretation, because of the displacement convexity of the functionals, is well suited to establish uniqueness and to study energy dissipation and convergence to equilibrium. A rigorous setting in which to make sense of the Otto calculus has been developed by Ambrosio, Gigli and Savar\'e \cite{AGS} for a suitable class of functionals. Convexity along geodesics in ${\cal W}_2$ also leads to new proofs of various  geometric and functional inequalities \cite{McCann}, \cite[Chapter 9]{Vil}. Finally, we mention that, when the space is not flat, qualitative properties of optimal transport can be quantified in terms of how bounds on the Ricci-Curbastro curvature affect the displacement convexity of certain specific functionals \cite[Part II]{Vil2}.

 The {\em tangent space} of $\mathcal P_2(\mR^N)$ at a probability measure $\mu$, denoted by $T_{\mu}\mathcal P_2(\mR^N)$ \cite{AGS} may be identified with the closure in  $L^2_{\mu}$ of the span of $\{\nabla\varphi:\varphi\in C^\infty_c\}$, where $C^\infty_c$ is the family of smooth functions with compact support. It is naturally equipped with the scalar product of $L^2_{\mu}$.

\section{The Fokker-Planck equation as a gradient flow on Wasserstein space}\label{FPGR}

Let us review the variational formulation of the Fokker-Planck equation as a gradient flow on Wasserstein space \cite{JKO,Vil,TGT}. Consider a physical system with  {\em phase space} $\mR^N$ and with {\em Hamiltonian} ${\cal H}: x\mapsto H(x)=E_x$. The thermodynamic states of the system are given by the family ${\cal P}(\mR^N)$ of probability distributions $P$ on $\mR^N$ admitting density $\rho$. On ${\cal P}(\mR^N)$, we define the {\em internal energy} as the expected value of the Energy {\em observable} in state $P$
\begin{equation}
U(H,\rho)=\E_P\{\mathcal H\}=\int_{\mR^N}H(x)\,\rho(x)dx=\langle H,\rho\rangle.
\end{equation}
Let us also introduce the (differential) {\em Gibbs entropy}
\begin{equation}
S(p)=-k\int_{\mR^N} \log \rho(x) \rho(x)dx,
\end{equation}
where $k$ is Boltzmann's constant. $S$ is strictly concave on ${\cal P}(\mR^N)$.  According to the Gibbsian postulate of classical statistical mechanics, the equilibrium
state of a microscopic system at constant absolute temperature T and with Hamiltonian function H is necessarily given by the Boltzmann distribution law with density
\begin{equation}\label{MB}\bar{\rho}(x)=Z^{-1}\exp\left[-\frac{H(x)}{kT}\right]
\end{equation}
where Z is the partition function\footnote{The letter $Z$ was chosen by Boltzmann to indicate ``zust\"{a}ndige Summe" (pertinent sum- here integral).}.  Let us introduce the {\em Free Energy} functional $F$ defined by
\begin{equation}\label{freeen}F(H,\rho,T):=U(H,\rho)-TS(\rho).
\end{equation}
Since $S$ is strictly concave on ${\cal S}$ and $U(E,\cdot)$ is linear, it follows that $F$ is strictly convex on the state space ${\cal P}(\mR^N)$. By Gibbs' variational principle, the Boltzmann distribution $\bar{\rho}$ is a minimum point of the free energy $F$ on ${\cal P}(\mR^N)$. Also notice that
\begin{eqnarray}\nonumber\D(\rho\|\bar{\rho})&=&\int_{\mR^N}\log\frac{\rho(x)}{\bar{\rho}(x)}\,\rho(x)dx\\&=&-\frac{1}{k}S(\rho)+\log Z+\frac{1}{kT}\int_{\mR^N}H(x)\rho(x)dx\nonumber\\&=&\frac{1}{kT}F(H,\rho,T)+ \log Z. \nonumber
\end{eqnarray}
Since $Z$ does not depend on $\rho$, we conclude that Gibb's principle is a trivial consequence of the fact that $\bar{\rho}$ minimizes $\D(\rho\|\bar{\rho})$ on ${\cal D}(\mR^N)$.

Consider now an absolutely continuous curve $\mu_t: [t_0,t_1]\rightarrow \mathcal W_2$. Then \cite[Chapter 8]{AGS}, there exist ``velocity field" $v_t\in L^2_{\mu_t}$ such that the following continuity equation holds on $(0,T)$
$$\frac{d}{dt}\mu_t+\nabla\cdot(v_t\mu_t)=0.$$
Suppose $d\mu_t=\rho_tdx$, so that the continuity equation
\begin{equation}\label{continuityeq}
\frac{\partial \rho}{\partial t}+\nabla\cdot(v\rho)=0
\end{equation}
holds. We want to study the free energy functional $F(H,\rho_t,T)$ or, equivalently, $\D(\rho_t\|\bar{\rho})$, along the flow $\{\rho_t; t_0\le t\le t_1\}$. Using (\ref{continuityeq}), we get
\begin{eqnarray}\nonumber
\frac{d}{dt}\D(\rho_t\|\bar{\rho})=\int_{\mR^N}\left[1+\log\rho_t+\frac{1}{kT}H(x)\right]\frac{\partial \rho_t}{\partial t}dx\\=-\int_{\mR^N}\left[1+\log\rho_t+\frac{1}{kT}H(x)\right]\nabla\cdot(v\rho_t)dx.\label{Decay}
\end{eqnarray}
Integrating by parts, if the boundary terms at infinity vanish, we get
\begin{eqnarray}\nonumber
\frac{d}{dt}\D(\rho_t\|\bar{\rho})=\int_{\mR^N}\nabla\left[\log\rho_t+\frac{1}{kT}H(x)\right]\cdot v \rho_tdx\\=\langle \nabla\log\rho_t+\frac{1}{kT}\nabla H(x), v\rangle_{L^2_{\rho_t}}.\nonumber
\end{eqnarray}
Thus, the Wasserstein gradient of $\D(\rho_t\|\bar{\rho})$ is
$$\nabla_{\mathcal W_2}\D(\rho_t\|\bar{\rho})=\nabla\log\rho_t+\frac{1}{kT}\nabla H(x).
$$
The corresponding gradient flow is
\begin{eqnarray}\nonumber
\frac{\partial\rho_t}{\partial t}=\nabla\cdot\left[\left(\nabla\log\rho_t+\frac{1}{kT}\nabla H(x)\right)\rho_t\right]\\=\nabla\cdot\left[\frac{1}{kT}\nabla H(x)\rho_t\right]+\Delta\rho_t.\label{FPflow}
\end{eqnarray}
But this is precisely the Fokker-Planck equation corresponding to the diffusion process
\begin{equation}\label{diffusionprocess}
dX_t=-\frac{1}{kT}\nabla H(X_t)dt+\sqrt{2}dW_t
\end{equation}
where $W$ is a standard $n$-dimensional Wiener process. The process (\ref{diffusionprocess}) has the Boltzmann distribution (\ref{MB}) as invariant density. Recall that \cite[p.220]{AGS} $F(H,\rho_t,T)$ or, equivalently, $\D(\rho_t\|\bar{\rho})$ are {\em displacement convex} and have therefore a unique minimizer.
\begin{remark}
{\em It seems worthwhile investigating to what extent the fundamental assumption of statistical mechanics that the variables with longer relaxation time form a {\em vector Markov process} having (\ref{MB}) as invariant density is equivalent to the requirement that the flow of one-time densities be a gradient flow in Wasserstein space for the free energy. }
\end{remark}

Let us finally plug the ``steepest descent" (\ref{FPflow}) into (\ref{Decay}). We get, after integrating by parts, the well known formula \cite{Gr}
\begin{eqnarray} \nonumber
\frac{d}{dt}\D(\rho_t\|\bar{\rho})=\int_{\mR^N}\left[1+\log\rho_t+\frac{1}{kT}H(x)\right]\frac{\partial \rho_t}{\partial t}dx\\\nonumber=\int_{\mR^N}\left[1+\log\rho_t+\frac{1}{kT}H(x)\right]\nabla\cdot\left[\frac{1}{kT}\nabla H(x)\rho_t+\nabla\rho_t\right]dx\\=-\int_{\mR^N}\|\nabla\log\left(\frac{\rho_t}{\bar{\rho}}\right)\|^2\rho_t dx.\label{FEdecay}
\end{eqnarray}
The last integral in (\ref{FEdecay}) is sometimes called the {\em relative Fisher information} of $\rho_t$ with respect to $\bar{\rho}$ \cite[p.278]{Vil}.

\section{Relative entropy as a functional on Wasserstein product spaces}\label{REWP}

Consider now two absolutely continuous curves $\mu_t: [t_0,t_1]\rightarrow \mathcal W_2$ and $\tilde{\mu}_t: [t_0,t_1]\rightarrow \mathcal W_2$ and their velocity fields $v_t\in L^2_{\mu_t}$ and $\tilde{v}_t\in L^2_{\tilde{\mu}_t}$. Then, on $(0,T)$
\begin{eqnarray}\label{conteq1}
\frac{d}{dt}\mu_t+\nabla\cdot(v_t\mu_t)=0,\\ \frac{d}{dt}\tilde{\mu}_t+\nabla\cdot(\tilde{v}_t\tilde{\mu}_t)=0.\label{conteq2}
\end{eqnarray}
Let us suppose that $d\mu_t=\rho_t(x)dx$ and $d\tilde{\mu}_t=\tilde{\rho}_t(x)dx$, for all $t\in [t_0,t_1]$. Then (\ref{conteq1})-(\ref{conteq2}) become
\begin{eqnarray}\label{Conteq1}
\frac{\partial \rho}{\partial t}+\nabla\cdot(v\rho)=0,\\ \frac{\partial\tilde{\rho}}{\partial t}+\nabla\cdot(\tilde{v}\tilde{\rho})=0,\label{Conteq2}
\end{eqnarray}
where the fields $v$ and $\tilde{v}$ satisfy
$$\int_{\mR^N}\|v(x,t)\|^2\rho_t(x)dx<\infty,\quad \int_{\mR^N}\|\tilde{v}(x,t)\|^2\tilde{\rho}_t(x)dx<\infty.
$$
The differentiability of the Wasserstein distance $W_2(\tilde{\rho}_t,\rho_t)$ has been studied \cite[Theorem 23.9]{Vil2}. Consider instead  the relative entropy functional
 on $\mathcal W_2\times \mathcal W_2$ 
 \begin{eqnarray}\nonumber \D(\tilde{\rho}_t\|\rho_t)=\int_{\mR^N}h(\tilde{\rho}_t,\rho_t)dx=\int_{\mR^N}\log\left(\frac{\tilde{\rho}_t}{\rho_t}\right)\tilde{\rho}_t dx,\\h(\tilde{\rho},\rho)=\log\left(\frac{\tilde{\rho}}{\rho}\right)\tilde{\rho}.\nonumber\end{eqnarray}
Relative entropy functionals $\D(\cdot\|\gamma)$, where $\gamma$ is a fixed probability measure (density), have been studied as geodesically convex functionals on $P_2(\mR^N)$, see \cite[Section 9.4]{AGS}.  Our study of the evolution of $\D(\tilde{\rho}_t\|\rho_t)$ is motivated by problems on a finite time interval such as the Schr\"{o}dinger bridge problem and stochastic control problems (Section \ref{stochastic control}) where it is important to evaluate relative entropy on {\em two} flows of marginals.

We  get
\begin{eqnarray}\nonumber
\frac{d}{dt}\D(\tilde{\rho}_t\|\rho_t)=\int_{\mR^N}\left[\frac{\partial h}{\partial \tilde{\rho}} \frac{\partial\tilde{\rho}}{\partial t}+\frac{\partial h}{\partial \rho} \frac{\partial \rho}{\partial t}\right]dx\\=\int_{\mR^N}\left[\left(1+\log\tilde{\rho}_t-\log\rho_t\right)\left(-\nabla\cdot(\tilde{v}\tilde{\rho}_t\right)\right.\nonumber\\\left.+\left(-\frac{\tilde{\rho}_t}{\rho_t}\right)\left(-\nabla\cdot(v\rho_t\right)\right]dx\label{entropyevolution}
\end{eqnarray}
After an integration by parts, assuming that the boundary terms at infinity vanish, we get
\begin{eqnarray}\nonumber
\frac{d}{dt}\D(\tilde{\rho}_t\|\rho_t)&=&\int_{\mR^N}\left[\nabla \log\left(\frac{\tilde{\rho}_t}{\rho_t}\right)\cdot \tilde{v}\tilde{\rho}_t-\nabla \frac{\tilde{\rho}_t}{\rho_t}\cdot v\rho_t\right]dx\\&=&\int_{\mR^N}\left[\left(\begin{matrix}\nabla \log\left(\frac{\tilde{\rho}_t}{\rho_t}\right)\\-\nabla \frac{\tilde{\rho}_t}{\rho_t}\end{matrix}\right)\cdot \left(\begin{matrix}\tilde{v}\tilde{\rho}_t\\v\rho_t\end{matrix}\right)\right]dx.\label{1}
\end{eqnarray}
Notice that the last expression looks like
$$\left\langle \left(\begin{matrix}\nabla \log\left(\frac{\tilde{\rho}_t}{\rho_t}\right)\\-\nabla \frac{\tilde{\rho}_t}{\rho_t}\end{matrix}\right), \left(\begin{matrix}\tilde{v}\\v\end{matrix}\right)\right\rangle_{L^2_{\tilde{\rho}_t}\times L^2_{\rho_t}}.
$$
Thus, we identify the gradient of the functional $\D(\tilde{\rho}\|\rho)$ on $\mathcal W_2\times \mathcal W_2$ as
\begin{equation}\label{gradWasser}
\left(\begin{matrix}\nabla^1_{\mathcal W_2} \D(\tilde{\rho}\|\rho)\\\nabla^2_{\mathcal W_2} \D(\tilde{\rho}\|\rho)\end{matrix}\right)=\left(\begin{matrix}\nabla \log\left(\frac{\tilde{\rho}}{\rho}\right)\\-\nabla \frac{\tilde{\rho}}{\rho_t}\end{matrix}\right).
\end{equation}
Let us now compute the gradient flow on  $\mathcal W_2\times \mathcal W_2$ corresponding to gradient (\ref{gradWasser}). We get
\begin{equation}\label{steepest}\frac{\partial}{\partial t}\left(\begin{matrix}\tilde{\rho}_t\\\rho_t\end{matrix}\right)-\nabla\cdot\left(\begin{matrix}\nabla\log\left(\frac{\tilde{\rho}_t}{\rho_t}\right)\tilde{\rho}_t\\-\nabla\left(\frac{\tilde{\rho}_t}{\rho_t}\right)\rho_t\end{matrix}\right)=0.
\end{equation}
Since 
$$J_1=\nabla\log\left(\frac{\tilde{\rho}_t}{\rho_t}\right)\tilde{\rho}_t=\nabla\left(\frac{\tilde{\rho}_t}{\rho_t}\right)\rho_t=-J_2,
$$
we observe the remarkable property that in the ``steepest descent" (\ref{steepest}) on the product Wasserstein space the ``fluxes" are {\em opposite} and, therefore, $\frac{\partial\tilde{\rho}}{\partial t}=-\frac{\partial\rho}{\partial t}$. If we plug the steepest descent (\ref{steepest}) into (\ref{entropyevolution}), we get what appears to be a new formula
\begin{eqnarray}\nonumber\frac{d}{dt}\D(\tilde{\rho}_t\|\rho_t)=\int_{\mR^N}\left[\left(1+\log\tilde{\rho}_t-\log\rho_t+\frac{\tilde{\rho}_t}{\rho_t}\right)\frac{\partial\tilde{\rho}}{\partial t}\right]dx\\\nonumber=-\int_{\mR^N}\left[\|\nabla\log\left(\frac{\tilde{\rho}_t}{\rho_t}\right)\|^2 \tilde{\rho}_t+\|\nabla\left(\frac{\tilde{\rho}_t}{\rho_t}\right)\|^2\rho_t\right]dx\\=-\int_{\mR^N}\left[\left(1+\frac{\tilde{\rho}_t}{\rho_t}\right)\|\nabla\log\left(\frac{\tilde{\rho}_t}{\rho_t}\right)\|^2 \tilde{\rho}_t\right]dx,\label{REFF}
\end{eqnarray}
which should be compared to (\ref{FEdecay}).

Let us return to equation (\ref{1}). By multiplying and dividing by $\tilde{\rho}_t$ in the last term of the middle expression, we get
\begin{equation}\label{PT2006}
\frac{d}{dt}\D(\tilde{\rho}_t\|\rho_t)=\int_{\mR^N}\left[\nabla \log\left(\frac{\tilde{\rho}_t}{\rho_t}\right)\cdot \left(\tilde{v}-v\right)\right]\tilde{\rho}_t dx\end{equation} 
which is precisely the expression obtained in \cite[Theorem III.1]{PT}.

\section{Elements of Nelson-F\"{o}llmer kinematics of finite-energy diffusion processes}\label{NF}

Let $(\Omega,\mathcal F, \mP)$ be a complete probability space. A stochastic process $\{\xi(t); t_0\le t\le t_1\}$ is called a {\em finite-energy diffusion} with constant diffusion coefficient $\sigma^2 I_N$ if the paths $\xi(\omega)$ belong to $C\left([t_0,t_1];\mR^N\right)$ ($N$-dimensional continuous functions) and
\begin{equation}\label{ito-forward}
\xi(t)-\xi(s)=\int_s^t\beta(\tau)d\tau+\sigma [W_+(t)-W_+(s)], \quad t_0\le s<t\le t_1,
\end{equation}
where $\beta(t)$ is at each time $t$ a measurable function of the past $\{\xi(\tau); t_0\le \tau\le t\}$ and $W$ is a standard $N$-dimensional Wiener process. Moreover, the drift $\beta$ satisfies the finite energy condition
$$\E\left\{\int_{t_0}^{t_1}\|\beta\|^2d\tau\right\}<\infty.
$$
In \cite{Foe}, F\"{o}llmer has shown that a finite-energy diffusion also admits a reverse-time Ito differential. Namely, there exists a measurable function $\gamma(t)$ of the future $\{\xi(\tau); t\le\tau\le t_1\}$ called {\em backward drift} and another Wiener process $W_-$ such that
 \begin{equation}\label{ito-backward}
\xi(t)-\xi(s)=\int_s^t\gamma(\tau)d\tau+\sigma [W_-(t)-W_-(s)], \quad t_0\le s<t\le t_1.
\end{equation}
Moreover, $\gamma$ satisfies
$$\E\left\{\int_{t_0}^{t_1}\|\gamma\|^2d\tau\right\}<\infty.
$$
Let us agree that $dt$ always indicate a strictly positive variable. For any function $f:[t_0,t_1]\rightarrow \mR$ let $d_+f(t)=f(t+dt)-f(t)$ be the {\em forward increment} at time $t$ and let $d_-f(t)=f(t)-f(t-dt)$ be the {\em backward increment} at time $t$. For a finite-energy diffusion, F\"{o}llmer has also shown in \cite{Foe} that forward and backward drifts may be obtained as Nelson's conditional derivatives \cite{N1}
\begin{eqnarray}\nonumber
\beta(t)&=&\lim_{dt\searrow 0}\E\left\{\frac{d_+\xi(t)}{dt}|\xi(\tau), t_0\le\tau\le t\right\},\\\gamma(t)&=&\lim_{dt\searrow 0}\E\left\{\frac{d_-\xi(t)}{dt}|\xi(\tau), t\le\tau\le t_1\right\},
\nonumber
\end{eqnarray}
the limits being taken in $L^2_N(\Omega,\mathcal F, \mP)$. It was finally shown in \cite{Foe} that the one-time probability density $\rho_t(\cdot)$ of $\xi(t)$ (which exists for every $t>t_0$) is absolutely continuous on $\mR^N$ and the following duality relation holds $\forall t>0$
\begin{equation}\label{duality}
\E\left\{\beta(t)-\gamma(t)|\xi(t)\right\}=\sigma^2\nabla\log\rho(\xi(t),t), \quad {\rm a.s.}.
\end{equation}
Let us introduce the fields
$$b_+(x,t)=\E\left\{\beta(t)|\xi(t)=x\right\},\quad b_-(x,t)=\quad \E\left\{\gamma(t)|\xi(t)=x\right\}.
$$
Then, Ito's rule for the forward and backward differential of $\xi$ imply that $\rho_t$ satisfies the two Fokker-Planck equations
\begin{eqnarray}\label{FP1}\frac{\partial \rho}{\partial t}+\nabla\cdot(b_+\rho)-\frac{\sigma^2}{2}\Delta\rho=0,\\\frac{\partial \rho}{\partial t}+\nabla\cdot(b_-\rho)+\frac{\sigma^2}{2}\Delta\rho=0.\label{FP2}
\end{eqnarray}
Following Nelson, let us introduce the {\em current} and {\em osmotic} drift of $\xi$ by
\begin{equation}\label{currosmot}v(t)=\frac{\beta(t)+\gamma(t)}{2}, \quad u(t)=\frac{\beta(t)-\gamma(t)}{2},
\end{equation}
respectively. Clearly $v$ is similar to the classical velocity, whereas $u$ is the velocity due to the noise which tends to zero when $\sigma^2$ tends to zero. Let us also introduce
$$v(x,t)=\E\left\{v(t)|\xi(t)=x\right\}=\frac{b_+(x,t)+b_-(x,t)}{2}.
$$
Then, combining (\ref{FP1}) and (\ref{FP2}), we get 
\begin{equation}\label{cont}
\frac{\partial \rho}{\partial t}+\nabla\cdot(v\rho)=0,
\end{equation}
which has the form of a continuity equation expressing conservation of mass.
When $\xi$ is {\em Markovian} with $\beta(t)=b_+(\xi(t),t)$ and $\gamma(t)=b_-(\xi(t),t)$, (\ref{duality}) reduces to Nelson's relation
\begin{equation}\label{Nelsonduality}b_+(x,t)-b_-(x,t)=\sigma^2\nabla\log\rho_t(x).
\end{equation}
Then (\ref{cont}) holds with
\begin{equation}\label{currvel}
v(x,t)=b_+(x,t)-\frac{\sigma^2}{2}\nabla\log\rho_t(x).
\end{equation}

\section{Relative entropy production for controlled evolution}\label{stochastic control}
Consider on $[t_0,t_1]$ a finite-energy Markov process taking values in $\mR^N$ with forward Ito differential
\begin{equation}\label{uncontrolled}d\xi=b_+(\xi(t),t)dt+\sigma dW_+.
\end{equation}
Let $\rho_t(x)$ be the probability density of $\xi(t)$. Consider also the feedback controlled process $\xi^u$ with forward differential
\begin{equation}\label{controlled}d\xi^u=b_+(\xi^u(t),t)dt+u(\xi^u(t),t)dt+\sigma dW_+.
\end{equation}
Here the control $u$ is adapted to the past and is such that $\xi^u$ is a finite-energy diffusion. Let $\rho_t^u(x)$ be the probability density of $\xi^u(t)$. We are interested in the evolution of $\D(\rho_t^u\|\rho_t)$. By (\ref{currvel})-(\ref{cont}), the densities satisfy
\begin{eqnarray}\nonumber
&&\frac{\partial \rho}{\partial t}+\nabla\cdot(v\rho)=0,\\\nonumber&&v(x,t)=b_+(x,t)-\frac{\sigma^2}{2}\nabla\log\rho_t(x)\\&&\frac{\partial\rho^u}{\partial t}+\nabla\cdot(v^u\rho^u)=0,\nonumber\\&&v^u(x,t)=b_+(x,t)+u(x,t)-\frac{\sigma^2}{2}\nabla\log\rho^u_t(x).\nonumber
\end{eqnarray}
By (\ref{PT2006}), we now get
\begin{eqnarray}\nonumber
\frac{d}{dt}\D(\rho^u_t\|\rho_t)=\int_{\mR^N}\left[\nabla \log\left(\frac{\rho^u_t}{\rho_t}\right)\cdot \left(v^u-v\right)\right]\rho^u_t dx\\=\int_{\mR^N}\left[\nabla \log\left(\frac{\rho^u_t}{\rho_t}\right)\cdot \left(u-\frac{\sigma^2}{2}\nabla \log\left(\frac{\rho^u_t}{\rho_t}\right)\right)\right]\rho^u_t dx.\label{PTcontr}
\end{eqnarray} 

Suppose now $\rho_t^u=\rho_t^0$ is also uncontrolled and differs from $\rho_t$ only because of the initial condition at $t=t_0$. Then (\ref{PTcontr}) gives the well known formula generalizing (\ref{FEdecay})
\begin{equation}\label{decay}
\frac{d}{dt}\D(\rho^0_t\|\rho_t)=-\frac{\sigma^2}{2}\int_{\mR^N}\left[\nabla \log\left(\frac{\rho^0_t}{\rho_t}\right)\cdot \nabla \log\left(\frac{\rho^0_t}{\rho_t}\right)\right]\rho^u_t dx
\end{equation} 
which shows that two solutions of the same Fokker-Plank equation tend to get closer.

Consider now the situation where $\xi(t)$ represents a ``prior" evolution on $[t_0,t_1]$ and  the controlled evolution $\xi^{u^*}=\xi^*$ is the solution of the Schr\"{o}dinger bridge problem for a pair of initial and final marginals $\rho_0$ and $\rho_1$ \cite{F2,W}. Then $$u^*(\xi^*,t)=\sigma^2\nabla\log\varphi(\xi^*,t)$$ and the differential of $\xi^*$ is given by
\begin{equation}\label{bridge}
d\xi^*=b_+(\xi^*(t),t)dt+\sigma^2\nabla\log\varphi(\xi^*(t),t)dt+\sigma dW_+
\end{equation}
where $\varphi$ is space-time harmonic for the prior evolution, namely it satisfies
\begin{equation}\label{spacetime}
\frac{\partial \varphi}{\partial t}+b_+\cdot\nabla\varphi+\frac{\sigma^2}{2}\Delta\varphi =0.
\end{equation}
Let $\rho^\varphi$ be the density of $\xi^*$. Let us first consider the special case of the Schr\"{o}dinger bridge problem where relative entropy on path space is minimized under the only constraint that the initial marginal density be $\rho_0\neq\rho_{t_0}$. Then, the optimal control $u^*$ is identically zero  and the evolution of the relative entropy is given by (\ref{decay}). Consider instead the case of the problem where only the final marginal density  $\rho_1\neq\rho_{t_1}$ is imposed. In such case,
$$\rho^\varphi_t(x)=\rho_t(x)\varphi(x,t).
$$
Then (\ref{PTcontr}) gives
\begin{equation}\label{bridge relentr}
\frac{d}{dt}\D(\rho^\varphi_t\|\rho_t)=\frac{\sigma^2}{2}\int_{\mR^N}\left[\nabla \log\varphi\cdot \nabla \log\varphi\right]\rho^\varphi_t dx.
\end{equation}
This shows that $\D(\rho^\varphi_t\|\rho_t)$ increases up to time $t=t_1$. It represents the intuitive fact that the bridge evolution has to be as close as possible to the prior but the final value of the relative entropy must be the positive quantity $\D(\rho_1\|\rho_{t_1})$. Thus, $\D(\rho^\varphi_t\|\rho_t)$ approaches this positive quantity from below. Result (\ref{bridge relentr}) may be viewed as a reverse-time H-theorem, as the bridge and the reference evolution have the same backward drift \cite{F2}.

\end{document}